\documentclass[12pt,reqno]{amsart}
\usepackage{amssymb,amsfonts,amsthm,amsmath}
\usepackage{cite}
\usepackage[shortlabels]{enumitem}
\usepackage[left=1 in,top=1 in,right=1 in, bottom=1 in]{geometry}
\usepackage{graphicx}
\usepackage{color}
\usepackage{mathtools}
\mathtoolsset{showonlyrefs}
\usepackage[pagebackref=false]{hyperref}
\usepackage{tcolorbox}
\def\d{{\rm d}}
\def\eqdef{\stackrel{\rm def}{=}}
\definecolor{darkred}{rgb}{.70,.12,.20}

\definecolor{darkgreen}{rgb}{.20,.52,.14}

\definecolor{byz}{rgb}{.44,.16,.39}

\numberwithin{equation}{section}

\newtheorem{theorem}{Theorem}[section]
\newtheorem{lemma}[theorem]{Lemma}

\newtheorem{corollary}[theorem]{Corollary}
\newtheorem{assumption}[theorem]{Assumption}

\newtheorem{proposition}[theorem]{Proposition}
\newtheorem{hypothesis}[theorem]{Hypothesis}

\theoremstyle{remark}
\newtheorem{remark}[theorem]{Remark}
\newtheorem{example}[theorem]{\bf{Example}}
\usepackage{todonotes}

\newcommand{\osc}{\mathop{\mathrm{osc}}}

\newcommand{\beq}{\begin{equation}}
\newcommand{\eeq}{\end{equation}}
\newcommand{\beqs}{\begin{equation*}}
\newcommand{\eeqs}{\end{equation*}}
\newcommand{\ba}{\begin{array}}
\newcommand{\ea}{\end{array}}
\newcommand{\beas}{\begin{eqnarray*}}
\newcommand{\eeas}{\end{eqnarray*}}
\newcommand{\bea}{\begin{eqnarray}}
\newcommand{\eea}{\end{eqnarray}}
\newcommand{\bal}{\begin{align}}
\newcommand{\eal}{\end{align}}

\newcommand{\bals}{\begin{align*}}
\newcommand{\eals}{\end{align*}}

\newcommand{\tnum}{\rm(\roman*)}
\newcommand{\rnum}{\rm(\alph*)}

\newcommand{\R}{\ensuremath{\mathbb R}}

\newcommand{\bds}{\begin{displaystyle}}
\newcommand{\eds}{\end{displaystyle}}

\newcommand{\remove}[1]{} %-- ON
\renewcommand{\remove}[1]{#1} % OFF

\definecolor{darkred}{rgb}{.70,.12,.20}

\definecolor{darkgreen}{rgb}{.20,.52,.14}

\begin{document}

\title{
A class of anisotropic diffusion-transport equations in non-divergence form
}

\author{Luan Hoang$^{1,*}$}
\address{$^1$Department of Mathematics and Statistics,
Texas Tech University,
1108 Memorial Circle, Lubbock, TX 79409--1042, U. S. A.}

\email{luan.hoang@ttu.edu}

\author{Akif Ibragimov$^{2}$}
\address{$^{2}$Department of Mathematics and Statistics,
Texas Tech University,
1108 Memorial Circle, Lubbock, TX 79409--1042, U. S. A.}
\email{akif.abraguimov@ttu.edu}

\thanks{$^*$Corresponding author.}

\date{\today}
%35Q30: NSE 
%76F02: Turbulence
%35C20: asymptotic expansion for solutions of PDE
%76D05  Navier-Stokes equations for incompressible viscous fluids

% 35A09   	Classical solutions to PDEs
% 35Q35   	PDEs in connection with fluid mechanics
% 35Q86   	PDEs in connection with geophysics
% 35B50   	Maximum principles in context of PDEs
% 35B40   	Asymptotic behavior of solutions to PDEs


\begin{abstract} 
We generalize Einstein's probabilistic method for the Brownian motion to study compressible fluids in porous media. The multi-dimensional case is considered with general probability distribution functions.
By relating the expected displacement per unit time with the velocity of the fluid, we derive an anisotropic diffusion equation in non-divergence form that contains a transport term. 
Under the Darcy law assumption, a corresponding nonlinear partial differential equations for the density function is obtained.
The classical solutions of this equation are studied, and the maximum and strong maximum principles are established. 
We also obtain exponential decay estimates for the solutions for all time, and particularly, their exponential convergence as time tends to infinity.
Our analysis uses some transformations of the Bernstein-Cole--Hopf type which are explicitly constructed even for very general equation of state.
Moreover, the Lemma of Growth in time is proved and utilized in order to achieve the above decaying estimates.
\end{abstract}

\maketitle 
\tableofcontents 

\pagestyle{myheadings}\markboth{\sc L.~Hoang and A.~Ibragimov}
{\sc Anisotropic diffusion-transport equations in non-divergence form}

%%%%%%%%%
\section{Introduction}\label{intro}
The purpose of this paper is two fold:
\begin{enumerate}[label=(\arabic*)]
    \item  Developing a new model for diffusion and transport processes of fluid flows in porous media using the Einstein paradigm for the Brownian motion \cite{Einstein1905}, and
    \item Analyzing that model rigorously to obtain certain stability results.
\end{enumerate}

Regarding the first goal, we recall that modeling the filtration in porous media 
is traditionally based on the following three components \cite{MuskatBook,BearBook,Scheidegger1974}: 

\begin{enumerate}[label=\rnum]
    \item  Equation of continuity (material balance/mass conservation).
    
    \item  Equation of motion which usually is the Darcy law or one of its generalizations.
    
    \item  Equation of state which describes the relation between the pressure and density.
\end{enumerate}
These lead to partial differential equations (PDE) of the parabolic type (linear, quasi-linear, degenerate, etc.) for the pressure or density function (see \cite{AronsonLec1986,CHIK1,Barletta2012}).
Thanks to (a), they all naturally  appear in the divergence form. 
They have been studied for so long with a vast literature, see, e.g., \cite{AronsonLec1986,VazquezPorousBook} for Darcy flows, \cite{ChadamQin,Payne1999a,Payne1999b, StraughanBook,ABHI1,HI1,HIKS1,CHIK1,HK2,CH2,CHK2,CHK4,HK3} for Forchheimer flows, and references therein. 

Although the three equations above are deterministic, they, in fact, can be subject to stochastic perturbations \cite{Rubin1991,WangLandau2001}.
Taking this stochastic point of view into account, we propose an alternative approach to the first component (a) -- the mass conservation -- by revisiting and utilizing Einstein's probabilistic material balance equation \cite{Einstein1905}.
More specifically, we apply Einstein's paradigm from \cite{Einstein1905} and consider 
fluid motion in porous media as events of random movements of the particles from point $x$ to  point $x+\zeta$ within the small time interval $\tau$, where $\zeta$ is a random displacement.
Generalizing the arguments in \cite{Einstein1905} to the multi-dimensional space, we arrive at the following PDE for the density function $\rho$
\beq \label{earlyeq}
\frac{\partial \rho}{\partial t}
=\sum_{i,j=1}^n a_{ij} \frac{\partial^2 \rho}{\partial x_i\partial x_j} + \frac{E}{\tau} \cdot \nabla \rho, 
\eeq 
where $a_{ij}(x,t)$ are the the diffusion coefficients, and 
$E(x,t)= \int \zeta \phi(x,t,\zeta)\d \zeta$ with $\phi(x,t,\zeta)$ denoting the probability distribution of these events, see more details in Section \ref{model} below.
Above, $E/\tau$ is a new transport term thanks to $\phi(x,t,\zeta)$ \emph{not} being assumed to be an even function in $\zeta$.
Note that $E$ is the expected value of the particle's displacement from time $t$ to $t+\tau$.
Therefore,  we postulate, see Hypothesis \ref{velassum}, that the average velocity $E/\tau$  is proportional to the velocity $v$ of the fluid, or, more generally,
\beq\label{as1} M_0 v =
\frac{1}{\tau} E ,
\eeq 
where $M_0$ is a matrix that guarantees that $v$ and  $E/\tau$ have a certain level of ``alignment", see \eqref{Evalign}.
This hypothesis \eqref{as1} connects the microscopic transport with the macroscopic one. It is essential to understanding and developing our  model.

After \eqref{earlyeq} and \eqref{as1}, 
we consider an anisotropic Darcy law for (b), and isentropic gas flows as well as slightly compressible fluid flows for (c).  They result in a quasi-linear parabolic equation of the second order in non-divergence form for $\rho$ that contains a quadratic term in $\nabla \rho$ and other nonlinearity in $\rho$, see \eqref{iseneq}--\eqref{drPeq} below.

We now turn to the second goal of the paper which is the mathematical analysis of the obtained models.
We  prove the maximum principle and strong maximum principle for the solutions. For the initial value problem with constant boundary data, we obtain exponential decay estimates for  the solution in the spatial $C^0$-norm.
Consequently, the solution exponentially converges in the $C^0$-norm, as time tends to infinity, to its constant boundary value. 
For the proofs, we explicitly construct \emph{separate} transformations of the Bernstein--Cole--Hopf type to convert the \emph{same} solution to a, whichever needed, sub- or super-solution of the corresponding truncated linear operator. Moreover, the Lemma of Growth in time is established  using the Landis method. It is then applied in consecutive time intervals to achieve the exponential decay estimates.

The paper is organized as follows.
In Section \ref{model} we derive the models in many steps. 
Firstly, by generalizing Einstein's probability method \cite{Einstein1905} to the multi-dimensional space,
 we derive a general diffusion equation \eqref{Deleq} of non-divergence form.
 \emph{Without} the assumption of the probability distribution function being even,
this equation contains the quotient $E(x,t)/\tau$ - the average displacement per unit time.
Secondly, by relating this quotient $E(x,t)/\tau$ with the velocity $v(x,t)$ of the fluid, we obtain equation \eqref{drPeq}.
The main assumption is Hypothesis \ref{velassum}, which generalizes the primary idea \eqref{mainEv}.
This hypothesis connects the microscopic concepts such as the particles movements with probabilities with the velocity which is a macroscopic feature of the fluid.
Thirdly, by using the ubiquitous Darcy law we find equation \eqref{eqgrav} for pressure $p$ and density $\rho$.
Lastly, the equation of state is used  to obtain the nonlinear PDE  \eqref{drPeq} for the density.
Special cases are equations \eqref{iseneq}, \eqref{idealeq} and \eqref{slighteq}.
Section \ref{maxmax} is devoted to studying equation \eqref{drPeq}, in its general form \eqref{maineq},  from the point of view of the maximum and strong maximum principles for the underlined nonlinear operator $L$, see  \eqref{Ldef}.
The maximum principle is proved in Theorem \ref{maxprin}, while the strong maximum principle is established in Theorem \ref{strmax}. For the latter theorem, we construct some transformations of the Bernstein--Cole--Hopf type in Lemma \ref{q-lin} to convert a solution of the operator $L$ to sub-solutions and super-solutions of the truncated linear operator $\mathcal L$. These transformations are explicit in terms of the equation of state.
In the last Section \ref{largesec}, we study the long-time behavior of the solution of the initial boundary value problem \eqref{nonIBVP}. The main tool is the Landis-typed Lemma of Growth for the linear operator $\widetilde L$, which is a general form of $\mathcal L$, in Lemma \ref{lemgrowth}. This leads to the estimates for the sub- and super- solutions of $\widetilde L$ in Proposition \ref{Hthm1}. With the aid  of the Bernstein--Cole--Hopf-typed transformations from Section \ref{maxmax}  to connect $L$, $\mathcal L$ and $\widetilde L$, we obtain the main results in Theorem \ref{Hthm2}. They consist of the exponential decay estimates  for all time, and, as a consequence,  the exponential convergence for the solutions as time tends to infinity.
In fact, the exponential rate can be independent of the initial data as shown in Corollary \ref{HCor}.
Applications to different types of fluid flows are given in Examples \ref{eg1} and \ref{eg2}.
It is worth mentioning that the Lemma of Growth manifests the stability of the original nonlinear problem which, in part, justifies our proposed framework for dynamics of fluid flows in porous media.

The authors are aware that the  models developed in this work are very different from the standard ones. Evidently, more data and experiments are needed to validate them. However, because the derivation is so straight forward, their ideas, methods and mathematical analysis are presented here with a hope that further discussions and developments will come along. They may eventually be useful in providing an alternative methodology to describe and understand complicated  diffusion and transport processes in porous media.
   
\section{Derivation of the models}\label{model}

\subsection*{Notation}
Throughout the paper, the spatial dimenion $n\ge 1$ is fixed. 
For a vector $x\in\R^n$, its Euclidean norm is denoted by $|x|$.
 Let $\mathcal M^{n\times n}$ denote the set of $n\times n$ matrices of real numbers, and $\mathcal M^{n\times n}_{{\rm sym}}$ denote the set of symmetric matrices in $\mathcal M^{n\times n}$.
For two matrices $A,B\in \mathcal M^{n\times n}$, their inner product $\langle A,B\rangle$ is the trace ${\rm Tr}(A^{\rm T}B)$. For a real-valued function $f(x)$ with $x=(x_1,\ldots,x_n)\in\R^n$, we denote by $D^2f$ the Hessian matrix of the second partial derivatives $(\partial^2f/\partial x_i\partial x_j)_{i,j=1,\ldots,n}$.

\subsection{General equations}
Let $\rho(x,t)$ be the density function of fluids at point $x\in \mathbb R^n$ and time $t\in\R$.
Let $\tau > 0$ be a small time interval as an input parameter at the time of observation.
Let $\zeta\in\R^n$ be the random displacement of the fluid particles.
Assume that the probability of the particles moving from location $x$ at time $t$ to location $x+\zeta$, for $\zeta=(\zeta_1,\zeta_2,\ldots,\zeta_n)\in\R^n$,  at time $t+\tau$ can be characterized by the 
probability distribution function $\phi(x,t,\zeta)\ge 0$ with  
  $$\int_{\R^n}\phi(x,t,\zeta)\d \zeta =1.$$  
The Einstein material balance \cite{Einstein1905} is
\begin{equation}\label{ein-matbal}
\rho(x,t+\tau)=\int_{\R^n} \rho(x+\zeta,t)\phi(x,t,\zeta)\d \zeta .
\end{equation}  

With $\tau$ small, we  approximate the time derivative of $\rho$ by
\begin{equation}\label{tau-assump}
 \frac{\partial\rho(x,t)}{\partial t}\approx \frac1{\tau}(\rho(x,t+\tau)- \rho(x,t))   .
\end{equation}

We will calculate $\rho(x,t+\tau)$ on the right-hand side of \eqref{tau-assump} by the material balance \eqref{ein-matbal}.
Assume the function $\zeta\mapsto \phi(x,t,\zeta)$ is supported in a small ball centered at the origin. By  the Taylor's expansion of the function $\zeta\mapsto \rho(x+\zeta,t)$, for small $|\zeta|$, up to the quadratic terms, we have the approximation
\begin{align*}
\rho(x+\zeta,t)\approx \rho(x,t)+\zeta \cdot \nabla \rho(x,t) +\frac12\sum_{i,j=1}^n \zeta_i\zeta_j \frac{\partial^2 \rho (x,t)}{\partial x_i\partial x_j} .
\end{align*}
Then using \eqref{ein-matbal}, we can approximate
\begin{align*}
    \rho(x,t+\tau)
    &=\int_{\R^n} \rho(x+\zeta,t)\phi(x,t,\zeta)\d \zeta 
    \approx \int_{\R^n}\rho(x,t)\phi(x,t,\zeta)\d \zeta \\
    &\quad + \int_{\R^n} \zeta \cdot \nabla \rho(x,t) \phi(x,t,\zeta) \d \zeta 
    +\frac12\sum_{i,j=1}^n   \int_{\R^n} \zeta_i\zeta_j \frac{\partial^2 \rho (x,t)}{\partial x_i\partial x_j}  \phi(x,t,\zeta)\d \zeta.
\end{align*}
Thus,
\beq\label{rapx}
    \rho(x,t+\tau)
    \approx 
 \rho (x,t)+ E(x,t) \cdot \nabla \rho(x,t) + \frac12\sum_{i,j=1}^n \bar a_{ij}(x,t) \frac{\partial^2 \rho(x,t)}{\partial x_i\partial x_j} ,
\eeq 
where the vector
\beq \label{Edef}
E(x,t)=\int_{\R^n} \phi(x,t,\zeta)\zeta \d \zeta,
\eeq 
and the coefficients
\beq\label{baraij}
    \bar a_{ij}(x,t)  =\int_{\R^n} \zeta_i\zeta_j \phi(x,t,\zeta) \d \zeta\text{ for } i,j=1,\ldots,n.
\eeq 

Combining \eqref{rapx} with \eqref{tau-assump} and replacing the approximation by equation, we obtain
\beq\label{Deleq0}
\frac{\partial \rho(x,t)}{\partial t}=\frac1{2\tau }\sum_{i,j=1}^n \bar a_{ij}(x,t) \frac{\partial^2 \rho(x,t)}{\partial x_i\partial x_j}  + \frac1{\tau }E(x,t) \cdot  \nabla \rho(x,t).
\eeq

Define the $n\times n$ matrices
\beq \label{Abar}
\bar A(x,t)=(\bar a_{ij}(x,t))_{i,j=1,\ldots,n}\text{ and } A(x,t)=(a_{ij}(x,t))_{i,j=1,\ldots,n}\eqdef \frac1{2\tau} \bar A(x,t)
\eeq 
Then equation \eqref{Deleq0} can be rewritten as
\beq\label{Deleq}
\frac{\partial \rho}{\partial t}=\sum_{i,j=1}^n a_{ij}\frac{\partial^2 \rho}{\partial x_i\partial x_j}  + \frac1{\tau }E \cdot  \nabla \rho
=\langle A,D^2\rho\rangle  + \frac1{\tau }E \cdot  \nabla \rho.
\eeq

\begin{remark}\label{Ermk}
The following remarks are in order.
\begin{enumerate}[label=\rnum]
    \item\label{rma} Thanks to \eqref{baraij}, the matrix $\bar A(x,t)$  is symmetric and, hence, so is $A(x,t)$.
    Also, for  $\xi=(\xi_1,...,\xi_n)\in\R^n$, one has 
$$\xi^{\rm T} \bar A(x,t) \xi =\sum_{i,j=1}^n \xi_i \bar a_{ij}(x,t) \xi_j =\int_{\R^n} |\xi\cdot \zeta|^2 \phi(x,t,\zeta) \d \zeta \ge 0.$$
Hence,  $\bar A(x,t)$ is positive semi-definite. Since $\tau>0$, then it is clear from \eqref{Abar} that the matrix $A(x,t)$ is also positive semi-definite.

\item \label{rmb} If $\zeta\mapsto \phi(x,t,\zeta)$ is an even function, then, by \eqref{Edef}, $E(x,t)=0$ and we have
\beq\label{Deleq1}
\frac{\partial \rho}{\partial t}=\sum_{i,j=1}^n a_{ij}\frac{\partial^2 \rho}{\partial x_i\partial x_j}.
\eeq

\item \label{rmc} Consider the case of mutually independent events with respect to  the coordinates of the displacement $\zeta$, that is, 
  \beqs %\label{indep} 
  \phi (x,t,\zeta) = \phi_1(x,t,\zeta_1)\cdots \phi_n(x,t,\zeta_n),\text{ for $\zeta=(\zeta_1,\zeta_2,\ldots,\zeta_n)$, }
  \eeqs 
with each $\phi_i(x,t,\zeta_i)$ being a probability distribution function in the variable $\zeta_i\in \R$,
 for $i=1,2,\ldots,n$.
 Then we have
    $$\bar a_{ij}=\begin{cases}
    \bar\sigma_i \bar\sigma_j,& \text{ for }i\ne j,\\
    \bar\sigma_{i,i}^2 , & \text{ for }i=j,
             \end{cases}
    $$
    where
\beqs
\bar\sigma_i(x,t)=
     \int_{\R} s \phi_i(x,t,s)\d s,\quad \bar\sigma_{i,i}(x,t)=\left (\int_{\R} s^2\phi_i(x,t,s) \d s\right)^{1/2}.
\eeqs 

\item \label{rmd} Assume, in addition to \ref{rmc}, that each function $\phi_i(x,t,s)$, for $1\le i\le n$,  is even in $s\in\R$.
Then each $\sigma_i=0$, and, hence,
$\bar A(x,t)$ is the diagonal matrix ${\rm diag}[\bar \sigma_{1,1}^2,\bar\sigma_{2,2}^2,\ldots,\bar\sigma_{n,n}^2]$.
Since each $\bar \sigma_{i,i}$ is positive, one has, in this case, that the matrix $\bar A(x,t)$ is positive definite.
Moreover, one finds that the function $\zeta\mapsto \phi(x,t,\zeta)$ is even, then, by part \ref{rmb},  equation \eqref{Deleq} becomes
\beq\label{Deleq2}
\frac{\partial \rho}{\partial t}=\sum_{i=1}^n \frac{\bar\sigma_{i,i}^2}{2\tau}\cdot \frac{\partial^2 \rho}{\partial x_i^2}.
\eeq
This is the multi-dimensional version of the equation obtained by Einstein in \cite{Einstein1905}.
\end{enumerate}
\end{remark}

We will focus on studying the general equation \eqref{Deleq}, not \eqref{Deleq1} nor \eqref{Deleq2}.

\subsection{Main assumption}\label{main asump}
Observe that $E(x,t)$ is the expected displacement from location $x$ between the time $t$ and $t+\tau$. Thus, $E(x,t)/\tau$ 
is the quotient
$$
\frac{\text{average displacement}}{\text{time of travel}}
$$
which can be seen as the average velocity.
Therefore, with small $\tau$, we can approximate this quotient $E(x,t)/\tau$ 
by the velocity $v(x,t)$  of the fluid, i.e.,
\beq\label{mainEv}
\frac{E(x,t)}{\tau} \approx v(x,t) .
\eeq
However, we will assume a much more general relation than \eqref{mainEv}.

\begin{hypothesis}\label{velassum}
There is a dimensionless  $n\times n$ matrix $M_0(x,t)$ such that 
\begin{equation}\label{expvhyp}
 M_0(x,t) v(x,t) = \frac{E(x,t)}{\tau} . 
\end{equation}
and 
\beq \label{weakM}
\xi^{\rm T} M_0(x,t)\xi\ge 0\text{ for all }\xi\in\R^n.
\eeq 
\end{hypothesis}

Condition \eqref{weakM} indicates that the velocity of the fluid $v(x,t)$ and the quotient $E(x,t)/\tau$   have some ``alignment", that is,
\beq \label{Evalign}
v(x,t)\cdot \frac{E(x,t)}{\tau} \ge 0.
\eeq 

Hypothesis \ref{velassum} is our fundamental assumption. It links the microscopic feature of the particles' movement in the media with the macroscopic property of the fluid flow -- the velocity of the fluid to be exact in this case.

Combining \eqref{Deleq} with \eqref{expvhyp} gives
\beq\label{drveq}
\frac{\partial \rho}{\partial t}=\langle A(x,t),D^2\rho\rangle  + (M_0(x,t) v(x,t))\cdot \nabla \rho.
\eeq

In this equation, the term $\langle A(x,t),D^2\rho\rangle$ represents the diffusion in the non-divergence form, and the term $(M_0(x,t) v(x,t))\cdot \nabla \rho$ represents the convection.

\subsection{Fluid motions in porous media}
Let $p(x,t)$ be the pressure of the fluid.
Assume the anisotropic Darcy's law \cite{DarcyBook,BearBook}, 
\beq \label{Darcy}
v=-\bar K(x,t)(\nabla p - \rho \vec g), 
\eeq 
where $\bar K(x,t)$ is an $n\times n$ matrix, 
and $\vec g$ is the gravitational acceleration for $n=1,2,3$, and can be any constant vector for $n\ge 4$.

Combining \eqref{drveq} with \eqref{Darcy} yields
\beq\label{eqgrav}
\frac{\partial \rho}{\partial t}=\langle A(x,t),D^2\rho\rangle -(K_0 (x,t)\nabla p)\cdot \nabla \rho + \rho B_0(x,t)\cdot \nabla \rho,
\eeq
where
\beqs %\label{KBdef}
 K_0(x,t)=M_0(x,t)\bar K(x,t),\quad B_0(x,t)=M_0(x,t) \bar K(x,t)\vec g.
\eeqs 

Next, we use equations of state to relate the pressure $p$ and density $\rho$ in \eqref{eqgrav}.

\medskip\noindent
\textit{Case of isentropic gas flows.} We have $p=c\rho^\gamma$ with a constant $c>0$ and the specific heat ratio $\gamma\ge 1$. Then \eqref{eqgrav} becomes
\beq\label{iseneq}
\frac{\partial \rho}{\partial t}=\langle A(x,t),D^2\rho\rangle - c\gamma \rho^{\gamma-1} (K_0(x,t)\nabla \rho)\cdot \nabla \rho +  \rho B_0(x,t)\cdot \nabla \rho.
\eeq
Particularly, for ideal gases, one has $\gamma=1$ and equation \eqref{iseneq} reads as
\beq\label{idealeq}
\frac{\partial \rho}{\partial t}=\langle A(x,t),D^2\rho\rangle - c(K_0(x,t)\nabla \rho)\cdot \nabla \rho +  \rho B_0(x,t)\cdot \nabla \rho.
\eeq

\medskip\noindent
\textit{Case of slightly compressible fluids.} We have
$$\frac1\rho \frac{d\rho}{dp}=\kappa, \text{ where $\kappa$ is the small, positive, constant compressibility.}$$
Noting that $\nabla \rho=\kappa\rho \nabla p$, we rewrite \eqref{eqgrav} as 
\beq\label{slighteq}
\frac{\partial \rho}{\partial t}=\langle A(x,t),D^2\rho\rangle  - \frac{1}{\kappa\rho} (K_0(x,t)\nabla \rho)\cdot \nabla \rho + \rho B_0(x,t) \cdot \nabla \rho.
\eeq

In general, assume the equation of state
\beq \label{plaw}
p=P_0(\rho), \text{ where $P_0$ is a known, increasing function.}
\eeq 
Then equation \eqref{eqgrav} becomes a PDE for $\rho$ which is
\beq\label{drPeq}
\frac{\partial \rho }{\partial t}=\langle A(x,t),D^2\rho\rangle  -P_0'(\rho)(K_0(x,t) \nabla \rho )\cdot \nabla \rho  + \rho  B_0(x,t)\cdot \nabla \rho .
\eeq

In the next two sections, we will focus entirely  on the mathematical aspect of equation \eqref{drPeq}.

\section{Maximum and strong maximum principles}\label{maxmax}

Let $U$ be a non-empty, open, bounded subset of  $\R^n$ with boundary $\Gamma$ and closure $\bar U$. Let $T>0$. 
Denote $U_T=U\times (0,T]$ and the parabolic boundary $\Gamma_T=\overline{U_T}\setminus U_T$, where $\overline{U_T}=\bar U\times[0,T]$ is the closure of $U_T$ (in $\R^{n+1}$).

Let $A:U_T\to \mathcal M^{n\times n}_{{\rm sym}}$ with $A(x,t)=(a_{ij}(x,t))_{i,j=1,\ldots,n}$, $K:U_T\to \mathcal M^{n\times n}$ and  $B:U_T\to \R^n$ be given functions.
Throughout this section \ref{maxmax}, we assume that there exists a constant $c_0>0$ such that
\beq\label{Aelip}
\xi^{\rm T} A(x,t)\xi \ge c_0|\xi|^2\text{ for all $(x,t)\in U_T$ and all $\xi\in \R^n$.}
\eeq 

\emph{Hereafter, $J$ is an interval  in $\R$ with non-empty interior, and $P$ is a function in $C^1(J,\R)$ with the derivative 
\beq \label{Pd}
P'\in C(J,[0,\infty)).
\eeq 
}

For fluid flows in porous media, $P$ is related to the equation of state  \eqref{plaw}.
However, we will consider general functions $P$.
It is clear from \eqref{Pd} that $P$ is an increasing function on $J$.

We denote by $C_{x,t}^{2,1}(U\times I)$, with any interval $I$ of $\R$, the class of continuous functions $u(x,t)$ with continuous partial derivatives $\partial u/\partial t$, $\partial u/\partial x_i$, $\partial^2 u/\partial x_i\partial x_j$ for $i,j=1.\ldots,n$.

We  study the nonlinear equation \eqref{drPeq} in the form
\beq \label{maineq}
\frac{\partial u}{\partial t}-\langle A(x,t),D^2u\rangle +   u B(x,t)\cdot \nabla u+ P'(u) (K(x,t)\nabla u)\cdot \nabla u=0.
\eeq 

For isentropic, non-ideal gas flows and equation \eqref{iseneq} with $\gamma>1$, we can use
\beq \label{isenchoice}
J=[0,\infty),\quad P(s)=s^\gamma,\quad K=cK_0,\quad B=-B_0.
\eeq 

For ideal gas flows and equation \eqref{idealeq}, we can still use \eqref{isenchoice} with $\gamma=1$ for physical $u=\rho$ (the density), but can also alter $J=\R$ in \eqref{isenchoice} for mathematical $u$. Thus,
\beq\label{idealchoice} 
J=[0,\infty) \text{ or } J=\R,\quad  P(s)=s,\quad K=cK_0,\quad B=-B_0.
\eeq 

For example, for slightly compressible fluids and equation \eqref{slighteq}, we can use
\beq \label{slighchoice}
J=(0,\infty),\quad P(s)=\ln s,\quad K=\kappa^{-1}K_0,\quad B=-B_0.
\eeq

We define the nonlinear operator $L$ associated with the left-hand side of \eqref{maineq} explicitly as
\beq \label{Ldef}
Lu=\frac{\partial u}{\partial t}-\langle A(x,t),D^2u\rangle +   u B(x,t)\cdot \nabla u+ P'(u) (K(x,t)\nabla u)\cdot \nabla u
\eeq 
for any function $u\in C_{x,t}^{2,1}(U_T)$ with the range $u(U_T)$ being a subset of $J$. 

Below, we study the maximum and strong maximum principles associated with this nonlinear operator $L$.

\subsection{Maximum principle}\label{maxsec}

\begin{theorem}[Maximum principle]\label{maxprin}
Assume 
\beq \label{class}
u\in C(\overline{U_T})\cap C_{x,t}^{2,1}(U_T)\text{ and $u(U_T)\subset J$.}
\eeq 
\begin{enumerate}[label=\tnum]
    \item     If  $Lu\le 0$ on $U_T$, then
    \beq\label{max}
        \max_{\overline{U_T}} u =\max_{\Gamma_T} u.
    \eeq

    \item  If  $Lu\ge 0$ on $U_T$, then
    \beq\label{min}
        \min_{\overline{U_T}} u =\min_{\Gamma_T} u.
    \eeq
\end{enumerate}
\end{theorem}
\begin{proof}
Define 
$\widetilde b(x,t)= u(x,t) B(x,t)+ P'(u(x,t))K(x,t)\nabla u(x,t)$
and the operator, for a function $v$,
\beqs
\widehat Lv= \frac{\partial v}{\partial t}-\langle A(x,t),D^2v\rangle +  \widetilde b(x,t)\cdot \nabla v.
\eeqs
Note that $\widehat Lu=Lu$.

(i) In this case we have $\widehat Lu\le 0$, hence, by the standard Maximum Principle for the linear operator $\widehat L$ and function $u$, we obtain \eqref{max}.

(ii) In this case we have $\widehat Lu\ge 0$, hence, by the standard Maximum Principle for the linear operator $\widehat L$ and function $u$, we obtain \eqref{min}.
\end{proof}

Let $S\subset \R^{n+1}$ and  $u$ be a bounded function on $S$. Denote
$$  \osc_S u =  \sup_{S} u-\inf_{S} u.$$

\begin{corollary}[Oscillation]\label{osccor}
Let the function $u$ be as in \eqref{class}.
  If $Lu=0$ on $U_T$, then 
    \beq\label{osc}
  \osc_{\overline{U_T}} u  = \osc_{\Gamma_T} u.
    \eeq
\end{corollary}
\begin{proof}
    Because $Lu=0$, we can apply both (i) and (ii) in Theorem \ref{maxprin}. Hence, we obtain \eqref{osc} from \eqref{max} and \eqref{min}.
\end{proof}

\subsection{Transformations of the Bernstein--Cole--Hopf type}\label{CHsec}

To remove the quadratic terms of the gradient, we introduce the following transformation of the Bernstein--Cole--Hopf type. 

For a given function $u$, we define an operator $\mathcal L$ as follows 
\beq \label{genL0}
 \mathcal Lw=\frac{\partial w}{\partial t}-\langle A(x,t),D^2w\rangle  +u(x,t) B(x,t)\cdot \nabla w.
\eeq  
Note that $\mathcal L$ is a linear operator in $w$ for each given function $u$.

\begin{lemma}\label{q-lin}
Let $u$ be a function such that 
\beq \label{class2}
u\in C_{x,t}^{2,1}(U_T)\text{ and $u(U_T)\subset J$,}
\eeq 
Define the linear operator $\mathcal L$ by \eqref{genL0}.
Given $s_0\in J$. For $\lambda\in\R$, $C>0$, $C'\in\R$, define
\beq \label{ftrans}
F_\lambda(s)=C\int_{s_0}^s e^{\lambda P(z)}dz +C'\text{ for }  s\in J.
\eeq 
\begin{enumerate}[label=\tnum]
    \item \label{Fprime} 
Then $F_\lambda\in C^2(J)$, $F_\lambda'>0$ and $\lambda  F_\lambda''\ge 0$ on $J$.
    
    \item\label{Fsub} 
Assume there is a  constant $c_1\ge 0$ such that 
\beq\label{cond1}
 \xi^{\rm T} K(x,t)\xi \ge -c_1|\xi|^2 \text{ for all $(x,t)\in U_T$, all $\xi\in\R^n$,} 
\eeq
If $Lu\le 0$ on $U_T$, then for any numbers $\lambda\ge c_1/c_0$, $C>0$, $C'\in\R$, the function $w=F_\lambda(u)$ satisfies
$\mathcal Lw\le 0$  on $U_T$.

 \item\label{Fsuper} 
Assume there is a constant $c_2\ge 0$ such that 
       \beq\label{cond2}
 \xi^{\rm T} K(x,t)\xi\le c_2|\xi|^2 \text{ for all $(x,t)\in U_T$, all $\xi\in\R^n$.}
   \eeq 
If $Lu\ge 0$ on $U_T$,  then for any numbers $\lambda\le  -c_2/c_0$, 
$C>0$, $C'\in\R$, the function $w=F_\lambda(u)$  satisfies
$\mathcal Lw\ge 0$  on $U_T$.
\end{enumerate}
\end{lemma}
\begin{proof}
\ref{Fprime} These facts clearly follow from \eqref{ftrans} and condition \eqref{Pd}.

For parts \ref{Fsub} and \ref{Fsuper}, we find  $w=F(u)$ for a function $F\in C^2(J)$ with 
\beq\label{Fder}
F'>0 \text{ on }J.
\eeq
    We have
    $$\frac{\partial w}{\partial x_i}=F'(u) \frac{\partial u}{\partial x_i},\quad 
    \frac{\partial^2 w}{\partial x_i\partial x_j}=F'(u)\frac{\partial^2 u}{\partial x_i\partial x_j}
    +F''(u)\frac{\partial u}{\partial x_i} \frac{\partial u}{\partial x_j}.$$
    Then
    \begin{align*}
 \mathcal Lw&=F'(u)\left[\frac{\partial u}{\partial t}-\sum_{i,j=1}^n a_{ij}\frac{\partial^2 u}{\partial x_{i}\partial x{_j}} +u B\cdot \nabla u \right]-F''(u)\sum_{i,j=1}^n a_{ij}\frac{\partial u}{\partial x_i} \frac{\partial u}{\partial x_j} \\
&=     F'(u)\mathcal Lu - F''(u) (\nabla u)^{\rm T} A(x,t)\nabla u.
    \end{align*}
 Note that 
$Lu = \mathcal Lu + P'(u)(\nabla u)^{\rm T} K(x,t)\nabla u $.
  Thus,
\beq\label{LLrel}
 \mathcal Lw=   F'(u)\{ Lu- P'(u)(\nabla u)^{\rm T} K(x,t)\nabla u\} - F''(u) (\nabla u)^{\rm T} A(x,t)\nabla u.
\eeq

\medskip
\ref{Fsub} Consider $Lu\le 0$ on $U_T$. Then it follows \eqref{LLrel} that
     \beqs
 \mathcal Lw\le   -P'(u) F'(u) (\nabla u)^{\rm T} K(x,t)\nabla u - F''(u) (\nabla u)^{\rm T} A(x,t)\nabla u.
    \eeqs 
 For $\mathcal Lw\le 0$ on $U_T$, we impose the condition
 \beq \label{presub}
   F''(u) \xi^{\rm T}  A(x,t)\xi\ge - P'(u) F'(u)    \xi^{\rm T} K(x,t) \xi \quad\forall (x,t)\in U_T,\forall \xi\in\R^n.
 \eeq 
We will find $F$ such that 
\beq \label{doubleF}
F''\ge 0\text{ on }J.
\eeq 
This property and \eqref{Aelip}, \eqref{cond1} imply, for all $(x,t)\in U_T$ and $\xi\in\R^n$, that
 \beqs 
   F''(u) \xi^{\rm T}  A(x,t)\xi\ge c_0F''(u)|\xi|^2, \quad  - P'(u) F'(u)    \xi^{\rm T} K(x,t) \xi\le c_1 P'(u) F'(u) |\xi|^2 .
 \eeqs 
Thanks to these inequalities and \eqref{presub}, a sufficient condition for $\mathcal Lw\le 0$ on $U_T$ is 
  \beqs 
  c_0F''(s)\ge c_1 P'(s)F'(s).
  \eeqs 
 For $\lambda\ge c_1/c_0$, the above inequality will follow from the equation    
    \beq\label{FFeq}
    F''(s)= \lambda P'(s)F'(s),
    \eeq
which yields a solution
\beq \label{Fex}
F'(s)= e^{\lambda \int P'(s) \d s }=Ce^{\lambda P(s)}.
\eeq 
Here, we choose $C>0$ so that condition \eqref{Fder} is met.
Then we select solution $F=F_\lambda$ as in \eqref{ftrans}. 
With \eqref{Fder} already satisfied,
equation \eqref{FFeq} and property \eqref{Pd} imply the second requirement  \eqref{doubleF}.

\medskip
\ref{Fsuper} Consider $Lu\ge 0$ on $U_T$. 
    Then one has from \eqref{LLrel} that
       \beqs 
 \mathcal Lw \ge  -P'(u) F'(u) (\nabla u)^{\rm T} K(x,t)\nabla u - F''(u) (\nabla u)^{\rm T} A(x,t)\nabla u.
    \eeqs
Thanks to this inequality, a sufficient condition for $\mathcal Lw\ge 0$ on $U_T$ is 
 \beq \label{presup}
F''(u) \xi^{\rm T} A(x,t)\xi\le - P'(u) F'(u)  \xi^{\rm T} K(x,t)\xi    \quad\forall (x,t)\in U_T,\forall \xi\in\R^n.
    \eeq 
In this case, we will find $F$ such that 
\beq \label{doubleF2}
F''\le 0\text{ on }J.
\eeq 
This and \eqref{Aelip} and \eqref{cond2} yield, for all $(x,t)\in U_T$ and $\xi\in\R^n$, that
 \beqs 
   F''(u) \xi^{\rm T}  A(x,t)\xi\le c_0 F''(u) |\xi|^2, \quad 
   - P'(u) F'(u)    \xi^{\rm T} K(x,t) \xi\ge -c_2 P'(u) F'(u) |\xi|^2 .
 \eeqs 
Using these inequalities and \eqref{presup}, a sufficient condition for $\mathcal Lw\ge 0$ on $U_T$ is
\beqs 
 c_0F''(s)\le -c_2 P'(s) F'(s).
\eeqs 
For $\lambda\le -c_2/c_0\le 0$, we, again,  solve equation \eqref{FFeq} instead. Same as in part  \ref{Fsub}, we choose solution $F=F_\lambda$ in \eqref{ftrans}.
Again, requirements \eqref{Fder} and \eqref{doubleF2} are met thanks to \eqref{Fex} and \eqref{FFeq}.
\end{proof}

Note that the function $F_\lambda$ in Lemma \ref{q-lin} is continuous and strictly increasing on $J$.

\begin{example}\label{fcase}
We have the following examples of fluid flows.
\begin{enumerate}[label=\rnum]
    \item\label{Fe1} \emph{Case of isentropic, non-ideal gas flows.} Use the choice in \eqref{isenchoice}.
Then we can choose
   \beq\label{Fl1} F_\lambda(s)=\int_0^s e^{\lambda z^\gamma}dz \text{ for } s\ge 0.
   \eeq 

    \item\label{Fe2}  \emph{Case of ideal gas.}
Use the choice in \eqref{idealchoice}. Choose $s_0=0$ in both cases of $J$.

For $\lambda=0$, we clearly can choose 
  \beq\label{Fl20} 
  F_\lambda(s)=s \text{ for }s\in J.
  \eeq 

For $\lambda\ne 0$, we can choose
  \beq\label{Fl21} 
  F_\lambda(s)=e^{\lambda s}\text{ for }s\in J.
  \eeq 

    \item\label{Fe3}  \emph{Case of slightly compressible fluids.} 
Use the choice in \eqref{slighchoice}.
In general, we have from \eqref{ftrans} that
    $$ F_\lambda(s)=C\int_{s_0}^s z^\lambda \d z+C'=\frac{C}{\lambda+1} (s^{\lambda+1}-s_0^{\lambda+1})+C'.$$ 
Therefore, we can choose 
\beq\label{Fl3} 
F_\lambda(s)=s^{\lambda+1} \text{ for } s>0.
\eeq 
\end{enumerate}  
\end{example}

\subsection{Strong maximum principle}\label{strongsec}
Assume additionally in this subsection \ref{strongsec} that $U$ is connected.

\begin{theorem}\label{strmax}
Assume both $A(x,t)$ and $B(x,t)$ are bounded on $U_T$, and  $K(x,t)$ satisfies 
condition \eqref{cond1} (respectively,  \eqref{cond2}.)
Suppose  $u$ is a function as in \eqref{class2},
 $u$ is bounded on $U_T$, and $Lu\leq 0$  (respectively, $Lu\geq 0$) on $U_T$.
Let 
$$M=\sup_{{U_T}} u(x,t)\text{ (respectively, } m=\inf_{{U_T}} u(x,t).)$$
Assume there is $(x_0,t_0)\in U_T$ such that 
  \beq \label{uMm} u(x_0, t_0)=M
   \quad (\text{respectively, } u(x_0, t_0)=m.)
   \eeq 
Then 
   \beq \label{uall}
   u(x,t)=M\quad  (\text{respectively, } u(x,t)=m) \text{ for all }(x,t)\in U_{t_0}=U\times(0,t_0].
   \eeq 
\end{theorem}
\begin{proof}
Let $\mathcal L w $ be defined by \eqref{genL0}. We rewrite $\mathcal L w $ as  
\beq \label{reL}
\mathcal L w = w_t-\langle A(x,t),D^2w\rangle  +\widetilde B(x,t)\cdot \nabla w,\text{ where } \widetilde B(x,t)=u(x,t)B(x,t). 
\eeq 

Because both $u(x,t)$ and  $B(x,t)$ are  bounded on $U_T$, we have $\widetilde B(x,t)$ is also bounded on $U_T$. Note that the operator $\mathcal L$ does not contain a term $c(x,t)w$.
Below, we use and refer to the  Strong Maximum Principle in the particular form  \cite[Chapter 3, Theorem 2.3]{LandisBook}, see also  \cite[Chapter 3, Theorem 2.4]{LandisBook}.

\medskip
\noindent\textit{Part 1.} Consider $Lu\le 0$ in $U_T$ and the corresponding conditions.  
Note from \eqref{uMm} that $u(x_0,t_0)=\max_{U_T} u(x,t)$. 

\emph{Case $c_1=0$.} We have $\mathcal Lu\le Lu\le 0$ on $U_T$ in this case. 
We can apply the Strong Maximum Principle to the operator $\mathcal L$, in the form of \eqref{reL},   and the function $u+|M|+1$ to have have $u=M$ on $U_{t_0}$.
 Thus, we obtain the first statement in \eqref{uall}.
 
\emph{Case $c_1>0$.} Let $\lambda_1=c_1/c_0$ and $w=F_{\lambda_1}(u)$ on $U_T$.
Thanks to Lemma \ref{q-lin}\ref{Fsub}, we have $\mathcal Lw\le 0$ on $U_T$,
and, thanks to the increase of $F_{\lambda_1}$,  
$$w(x_0,t_0) =F_{\lambda_1}(M)=\max_{U_T} w.$$ 
We apply the Strong Maximum Principle to the operator $\mathcal L$  and the function $w+|F_{\lambda_1}(M)|+1$ to have have $w=F_{\lambda_1}(M)$ on $U_{t_0}$.
Thus, $u=F_{\lambda_1}^{-1}(w)=M$ on $U_{t_0}$.

\medskip
\noindent\textit{Part 2.}  Consider $Lu\ge 0$  in $U_T$ and the corresponding conditions.
Note that $u(x_0,t_0)=\min_{U_T} u(x,t)$.

\emph{Case $c_2=0$.} We have $\mathcal Lu\ge Lu\ge 0$ on $U_T$.  
We apply the Strong Maximum Principle to the operator $\mathcal L$ and the function $u-|m|-1$ to have have $u=m$ on $U_{t_0}$.
 Thus, we obtain the second statement in \eqref{uall}.

 \emph{Case $c_2>0$.} Let $\lambda_2=-c_2/c_0$ and $w=F_{\lambda_2}(u)$ on ${U_T}$.
The proof of the second statement in \eqref{uall} is similar to Part 1, Case $c_1>0$ above with the use of Lemma \ref{q-lin}\ref{Fsuper} instead of Lemma \ref{q-lin}\ref{Fsub}, and the Strong Maximum Principle applied to the operator $\mathcal L$ and the function $u-|m|-1$. We omit the details.  
\end{proof}

\section{Initial boundary value problem}\label{largesec}

Let $U$ and $\Gamma$ be as in Section \ref{maxmax}.
We fix a point $x_0\not\in \bar U$ and set
\beq
r_0=\min\{|x-x_0|:x\in\bar U\},\quad R=\max\{|x-x_0|:x\in\bar U\}.
\eeq
Then $R>r_0>0$.

\subsection{Results for linear operators}

As seen in Section \ref{maxmax}, the study of the nonlinear problem can be reduced to that of some linear operator. 
Therefore, we first establish some essential results for a general  linear case. 

Given $T>0$, let $U_T$ and $\Gamma_T$ be as in Section \ref{maxmax}.

\begin{assumption}\label{firstA}
Let $A:U_T\to \mathcal M^{n\times n}_{{\rm sym}}$ and $b:U_T\to \R^n$ be such that
\begin{enumerate}[label=\tnum]
    \item $A$ satisfies \eqref{Aelip} for some constant $c_0>0$, and 
    \item there are constants $M_1>0$ and $M_2\ge 0$ such that
\beq\label{TAB}
{\rm Tr}(A(x,t))\le M_1,\quad |b(x,t)|\le M_2 \text{ for all }(x,t)\in U_T.
\eeq
\end{enumerate}
\end{assumption}

Define the linear operator $\widetilde L $ by
\beq \label{Ltil}
\widetilde L w = w_t-\langle A(x,t),D^2w\rangle  +b(x,t)\cdot \nabla w \text{ for $w\in C_{x,t}^{2,1}(U_T)$.} 
\eeq 

\begin{lemma}[Lemma of Growth]\label{lemgrowth}
Under Assumption \ref{firstA}, set 
 \beq\label{sTe}
 \beta=\frac{1}{4c_0}\max\left \{ 2(M_1+M_2R),\frac{R^2}{T}\right\}, \quad T_*=\frac{R^2}{4c_0 \beta},\quad 
 \eta_*=1-(r_0/R)^{2\beta}.
 \eeq
 If  $w\in C(\overline{U_T})\cap C_{x,t}^{2,1}(U_T)$ satisfies $\widetilde Lw\le 0$ on $U_T$ and $w\le 0$ on $\Gamma\times[0,T]$, then one has
 \beq \label{Tgrow}
\max\{0,\max_{x\in\bar U} w(x,T_*)\}\le \eta_*\max\{0,\max_{x\in\bar U} w(x,0)\}.
\eeq 
\end{lemma}
\begin{proof}
We follow \cite[Chapter 3, Lemma 6.1]{LandisBook} and also \cite[Lemma IV.3]{HIK2}.
We ignore the values in \eqref{sTe} momentarily.

\medskip\noindent\textit{Step 1.} 
Let $\varphi\in C(\bar U)\cap C^2(U)$ such that 
\beq \label{genphi}
0<d_0\le \varphi\le d_1\text{ on $\bar U$, and }
|\nabla \varphi|\le d_2,\  
\varphi\le c_0|\nabla \varphi|^2,\ 
|\langle A,D^2\varphi \rangle|\le d_3  \text{ on $U$,}
\eeq
for some positive numbers $d_0$, $d_1$, $d_2$, $d_3$. 
A specific function $\varphi$ will be constructed in Step 3 below.
Define a barrier function $W$ on $\bar U\times \R$ by
\beqs
W(x,t)=\begin{cases}
    t^{-\beta}e^{-\frac{\varphi(x)}{t}}& \text{ if } (x,t)\in \bar U\times(0,\infty),\\
    0,&\text{ if } (x,t)\in \bar U\times(-\infty,0],
\end{cases}
\eeqs
where $\beta$ is a positive number.
Thanks to the lower bound of $\varphi(x)$ in \eqref{genphi}, we have $W\in C(\bar U\times \R)$.
With the calculations, for $t>0$,
\beqs
W_t=-\frac{\beta}{t}W+\frac{\varphi}{t^2}W,\quad W_{x_i}=-\frac{\varphi_{x_i}}{t}W,\quad 
W_{x_ix_j}=-\frac{\varphi_{x_ix_j}}{t}W+\frac{\varphi_{x_i}\varphi_{x_j}}{t^2}W,
\eeqs
we have, on $U_T$,
\beqs
\widetilde L W=\frac{W}{t^2}\left\{ t(-\beta+ \langle A,D^2\varphi \rangle - b\cdot \nabla \varphi) +\varphi   -(A\nabla \varphi)\cdot\nabla \varphi \right\}.
\eeqs

We aim at having $\widetilde LW\le 0$ on $U\times (0,\infty)$, thus impose the conditions 
\beq\label{spineq}
\varphi  \le (A\nabla \varphi)\cdot\nabla \varphi
\text{ and } 
\beta\ge  \langle A,D^2\varphi \rangle - b\cdot \nabla \varphi 
 \text{ on } U_T.
\eeq
Thanks to \eqref{Aelip}, a sufficient condition for the first condition in \eqref{spineq} is
$\varphi\le c_0|\nabla \varphi|^2$ on $U_T$ which, in fact, is met by our assumption \eqref{genphi}.
A sufficient condition for the second condition in \eqref{spineq} is
\beq \label{sval}
\beta\ge |\langle A,D^2\varphi \rangle| + |b| |\nabla \varphi| \text{ on } U_T.
\eeq 
Based on \eqref{TAB} and \eqref{genphi}, we will choose 
\beq \label{sval0}
\beta\ge d_3 + M_2 d_2
\eeq 
in order to satisfy \eqref{sval}. In summary, for $\beta$ in \eqref{sval0}, we have $\widetilde L W\le 0$ on $U\times(0,\infty)$. 

\medskip\noindent\textit{Step 2.} 
Let $\beta$ satisfy \eqref{sval0}.
Set $M=\max\{0,\max_{\bar U} w(x,0)\}$
and define $\widetilde W=M(1-\eta W)$ on $\bar U\times \R$, with $\eta=(d_0 e/\beta)^\beta>0$.
Then $\widetilde L \widetilde W\ge 0$ on $U\times(0,\infty)$. 
On $\bar U\times\{0\}$, one has $W(x,0)=0$ and, hence, 
\beq\label{bdv1}
\widetilde W(x,0)=M\ge w(x,0)\text{ for all $x\in \bar U$.}
\eeq
Observe, for $t>0$ that 
\beqs
\widetilde W(x,t)=M\left (1-\eta t^{-\beta} e^{-\varphi(x)/t}\right)\ge M\left(1-\eta t^{-\beta}e^{-d_0/t}\right).
\eeqs
Elementary calculations show that the function $h_0(t)=t^{-\beta}e^{-d_0/t}$ on $(0,\infty)$ attains  the maximum at $t_0=d_0/\beta$ with the value $h_0(t_0)=\eta^{-1}$.
Thus, one has, on $\bar U\times(0,\infty)$, 
\beqs
\widetilde W(x,t)\ge M\left(1-\eta h_0(t_0)\right)=0.
\eeqs
In particular,
\beq\label{bdv2}
\widetilde W\ge w \text{ on }\Gamma\times(0,\infty).
\eeq

We impose another condition  
\beq\label{newsT}
\beta\ge d_1/T\text{ and set }T_*=d_1/\beta.
\eeq 
Then $T_*\le T$, and one has $\widetilde L(w-\widetilde W)\le 0$ on $U\times(0,T_*]$, and, thanks to \eqref{bdv1}, \eqref{bdv2}, $w-\widetilde W \le 0$ on the parabolic boundary of $U\times(0,T_*]$. Applying the Maximum Principle to the operator $\widetilde L$  and function $(w-\widetilde W)$ on the set $\bar U\times [0,T_*]$ yields
\beq\label{Ww}
w\le \widetilde W \text{ on } \bar U\times [0,T_*].
\eeq
Note that
\beq\label{We}
\widetilde W(x,t)\le M\left(1-\eta t^{-\beta} e^{-d_1/t}\right).
\eeq

When $t=T_*$, it follows \eqref{Ww} and \eqref{We}, for all $x\in \bar U$, that 
\beqs
w(x,T_*)\le \widetilde W(x,T_*)
\le M\left[1-\left(\frac{d_0e}{\beta}\right )^\beta \left(\frac{d_1}{\beta}\right)^{-\beta} e^{-d_1(\beta/d_1)}\right]
=M\left[1-\left(\frac{d_0}{d_1}\right)^{\beta}\right]=\eta_* M,
\eeqs
where 
\beq\label{neweta}
\eta_*=1-(d_0/d_1)^{\beta}\in(0,1).
\eeq
Thus, we obtain inequality \eqref{Tgrow} with $T_*$, $\eta_*$ together with $\beta$ as in \eqref{sval0}, \eqref{newsT}, \eqref{neweta}.

\medskip\noindent\textit{Step 3.} 
Specifically, we choose the function $\varphi(x)=\mu |x-x_0|^2$, for a number $\mu>0$ chosen later. Then
$\mu r_0^2\le \varphi\le \mu R^2$ on $\bar U$, hence, we choose
$d_0=\mu r_0^2$ and $d_1=\mu R^2$ in \eqref{genphi}.

For the second condition in \eqref{genphi}, since $|\nabla \varphi|=2\mu|x-x_0|\le 2\mu R$ we choose $d_2=\mu R$.

The third condition in \eqref{genphi} becomes $\mu |x-x_0|^2 \le 4 c_0 \mu^2 |x-x_0|^2$, thus,  we select
\beqs %\label{muval}
\mu=\frac1{4c_0}.
\eeqs

For last condition in \eqref{genphi}, we note that $\langle A,D^2\varphi \rangle= 2\mu {\rm Tr}(A(x,t))$, and, thus, choose $d_3=2\mu M_1$.

With the above values of $\mu$, $d_0$, $d_1$, $d_2$, $d_3$, the relations  \eqref{sval0} and \eqref{newsT} become
\beq\label{slast}
\beta\ge 2\mu (M_1+M_2R)=\frac{2(M_1+M_2R)}{4c_0}, \quad \beta\ge \frac{\mu R^2}{T}=\frac{R^2}{4c_0 T},\quad T_*=\frac{\mu R^2}{\beta}=\frac{R^2}{4c_0 \beta}.
\eeq

For $\beta$ chosen in \eqref{sTe}, it satisfies the first two conditions in \eqref{slast}. 
Also, $T_*$ in \eqref{slast} is exactly the one given in \eqref{sTe}. Moreover, we have from \eqref{neweta} that 
$$\eta_*=1-(\mu r_0^2/ (\mu R^2))^\beta=1-(r_0/R)^{2\beta}$$
which is the same number in \eqref{sTe}. The proof is complete.
\end{proof}

The key improvement in estimate \eqref{Tgrow} compared with the Maximum Principle  is the factor $\eta_*$ belonging to $(0,1)$.
From the above Lemma of Growth, we derive more specific estimates for sub/super-solutions and the solutions themselves for all time.
The emphasis is the decaying ones for large time, even though some ``optimal" estimates for small time are also obtained.

\begin{assumption}\label{secondA}
Let  $A:U\times(0,\infty)\to \mathcal M^{n\times n}_{{\rm sym}}$ and $b:U\times(0,\infty)\to \R^n$ satisfy
\begin{enumerate}[label=\tnum]
    \item there is a positive constant $c_0$ such that 
\beq\label{Aelip2}
\xi^{\rm T} A(x,t)\xi \ge c_0|\xi|^2\text{ for all $(x,t)\in U\times(0,\infty)$ and all $\xi\in \R^n$,}
\eeq 
and
    \item   $A(x,t)$ and  $b(x,t)$ are bounded on $U\times (0,\infty)$.
\end{enumerate}
\end{assumption}

Under Assumption \ref{secondA}, define the linear operator $\widetilde L$ by \eqref{Ltil} for $w\in C_{x,t}^{2,1}(U\times(0,\infty))$.
Thanks to condition (ii) in Assumption \ref{secondA}, there are constants $M_1>0$ and $M_2>0$ such that
\beq\label{TAB2}
{\rm Tr}(A(x,t))\le M_1,\quad |b(x,t)|\le M_2 \text{ for all }(x,t)\in U\times(0,\infty).
\eeq

\begin{proposition}\label{Hthm1}
Let Assumption \ref{secondA} hold and positive numbers $M_1$, $M_2$ be as in \eqref{TAB2}.
Set
\beq\label{sten}
\begin{aligned}
 \beta&=\frac1{2c_0}(M_1+M_2R), \
 T_*=\frac{R^2}{4c_0 \beta},\ 
 \eta_*=1-(r_0/R)^{2\beta},\\
 \nu&=T_*^{-1}\ln(1/\eta_*),\
 \nu_0=\frac{R^2}{2c_0}\ln (R/r_0).
\end{aligned} 
\eeq
Let $w\in C(\bar U\times[0,\infty))\cap C_{x,t}^{2,1}(U\times(0,\infty))$.
\begin{enumerate}[label=\tnum]
    \item If  $\widetilde Lw\le 0$ on $U\times (0,\infty)$ and $w\le 0$  on $\Gamma\times(0,\infty)$, then
     \begin{align}\label{upper1} 
\max_{x\in \bar U} w(x,t)&\le (1-e^{-\nu_0/ t})\max\{0,\max_{x\in \bar U} w(x,0)\} &&\text{ for } 0< t\le T_*,\\
\label{upper2}
    \max_{x\in \bar U} w(x,t) &\le \eta_*^{-1} e^{-\nu t }\max\{0,\max_{x\in \bar U} w(x,0)\} &&\text{ for }t\ge 0,
    \end{align}
    and, consequently,
    \begin{equation}\label{basic1}
    \limsup_{t\to \infty} \max_{x\in \bar U} w(x,t) \le 0.
    \end{equation}

    \item If  $\widetilde Lw\ge 0$  on $U\times (0,\infty)$ and $w\ge 0$  on $\Gamma\times(0,\infty)$, then
         \begin{align}\label{lower1} 
\min_{x\in \bar U} w(x,t)&\ge (1-e^{-\nu_0/ t})\min\{0,\min_{x\in \bar U} w(x,0)\} &&\text{ for } 0< t\le T_*,\\
\label{lower2}
    \min_{x\in \bar U} w(x,t) &\ge \eta_*^{-1} e^{-\nu t }\min\{0,\min_{x\in \bar U} w(x,0)\} &&\text{ for }t\ge 0,
    \end{align}
    and, consequently,
    \begin{equation}\label{basic2}
    \liminf_{t\to \infty} \min_{x\in \bar U} w(x,t) \ge 0.
    \end{equation}

    \item If  $\widetilde Lw= 0$  on $U\times (0,\infty)$ and $w=0$  on $\Gamma\times(0,\infty)$, then
    \begin{align}\label{abs1}
    \max_{x\in \bar U} |w(x,t)| &\le (1-e^{-\nu_0/ t})\max_{x\in \bar U}|w(x,0)| &&\text{ for } 0< t\le T_*,\\
    \max_{x\in \bar U} |w(x,t)| &\le \eta_*^{-1} e^{-\nu t }\max_{x\in \bar U}|w(x,0)|&&\text{ for }t\ge 0,\label{abs2}
    \end{align}
        and, consequently,
    \begin{equation}\label{basic3}
    \lim_{t\to \infty} \max_{x\in \bar U} |w(x,t)| = 0.
\end{equation}
  \end{enumerate}  
\end{proposition}
\begin{proof} 
For any integer $k\ge 0$, let $$T_k=kT_*,\quad  J_k=\max\{0,\max_{x\in \bar U} w(x,T_k)\}\ge 0.$$

\medskip
(i) In this case, $\widetilde Lw\le 0$ on $U\times (0,\infty)$ and $w\le 0$  on $\Gamma\times[0,\infty)$.

\medskip 
We prove \eqref{upper1} first. Letting $t\in(0,T_*]$, we apply Lemma \ref{lemgrowth} to $T=t$. Using $\beta'$, $T_*'$, $\eta_*'$ to denote $\beta$, $T_*$, $\eta_*$ in \eqref{sTe}, and noticing 
$$t\le  T_*=\frac{R^2}{4c_0 \beta}=\frac{R^2}{2(M_1+M_2R)},$$
we have 
\begin{align*} 
\beta'&=\frac{1}{4c_0}\max\left \{2(M_1+2M_2R),\frac{R^2}{t}\right\}=\frac{R^2}{4c_0 t},\\
T_*'&=R^2/(4c_0s')=t,\text{ and } 
\eta_*'=1-(r_0/R)^{R^2/(2c_0t)}=1-e^{-\nu_0/ t}.
\end{align*}
We obtain from \eqref{Tgrow} that 
 \beqs %\label{up11}
\max_{\bar U} w(x,t)\le \eta_*'\max\{0,\max_{\bar U} w(x,0)\}=(1-e^{-\nu_0/ t})\max\{0,\max_{\bar U} w(x,0)\},
 \eeqs 
which proves \eqref{upper1}.

\medskip
We prove \eqref{upper2} next.  For $k\ge 1$, we apply Lemma \ref{lemgrowth} to the cylinder $\bar U\times [T_{k-1},T_k]$  to have
$J_k\le \eta_* J_{k-1}$.
Iterating this inequality in $k$ gives
\beq\label{Jineq}
J_k\le \eta_*^k J_0\text{ for any }k\ge 0.
\eeq 
For each $t>0$, let $k\ge 0$ such that $t\in(T_k,T_{k+1}]$. 
Note  $k+1\ge t/T_*$.
By the Maximum Principle, the fact that $w\le 0$ on $\Gamma\times(T_{k-1},T_k]$, and inequality \eqref{Jineq}, one has
\beqs
w(x,t)\le J_k\le \eta_*^k J_0=\eta_*^{-1}\eta_*^{k+1} J_0\le \eta_*^{-1} \eta_*^{t/T_*}J_0=\eta_*^{-1} e^{-t T_*^{-1}\ln(1/\eta_*)}J_0.
\eeqs
Thus,
 \beqs %\label{up22}
    \max_{x\in \bar U} w(x,t)\le \eta_*^{-1} e^{-\nu t }J_0 \text{ for any }t\ge 0,
 \eeqs  
which proves \eqref{upper2}.
Passing $t\to\infty$ in \eqref{upper2} yields \eqref{basic1}.

\medskip
(ii) In this case, we can apply the results in part (i) by replacing $w$ with $-w$.
Then it follows \eqref{upper1} and \eqref{upper2} that
\begin{equation}\label{up1}
    \max_{x\in \bar U} (-w(x,t))\le (1-e^{-\nu_0/ t })\max\{0,\max_{\bar U} (-w(x,0))\} \text{ for }0< t\le T_*,
    \end{equation}
\begin{equation}\label{up2}
    \max_{x\in \bar U} (-w(x,t))\le \eta_*^{-1} e^{-\nu t }\max\{0,\max_{\bar U} (-w(x,0))\} \text{ for }t\ge 0,
    \end{equation}
  which imply \eqref{lower1} and \eqref{lower2}, respectively. Passing $t\to\infty$ in \eqref{lower2} yields \eqref{basic2}.

\medskip
(iii) Since $\widetilde Lw=0$  on $U\times (0,\infty)$ and $w=0$ on $\Gamma\times(0,\infty)$, we can apply the results in both parts (i) and (ii) above.
Observe that
\beq\label{absw} 
|w(x,t)|=\max\{ w(x,t),-w(x,t)\}\le  \max\{ \max_{x\in \bar U} w(x,t),\max_{x\in \bar U} (-w(x,t))\}.
\eeq 

For $0<t\le T_*$, combining \eqref{absw} with \eqref{upper1} and \eqref{up1}, and using the fact
\beq\label{clearini}
\max_{\bar U} w(x,0),\  \max_{\bar U} (-w(x,0))\le \max_{\bar U} |w(x,0)|,
\eeq
we obtain \eqref{abs1}.

For $t\ge 0$, combining \eqref{absw} with \eqref{upper2}, \eqref{up2} and \eqref{clearini} gives \eqref{abs2}.
Finally, \eqref{basic3} follows \eqref{abs2}.
\end{proof}

\begin{remark}
Note that the estimates \eqref{upper1}, \eqref{lower1} and \eqref{basic1} for small time  are more optimal, as $t\to 0^+$, than their counterparts for large time \eqref{upper2}, \eqref{lower2} and \eqref{basic3}. It is because the factors in front of the initial data converge to $1$, as $t\to 0^+$, instead of $\eta_*^{-1}>1$.
\end{remark}

\subsection{Results for the nonlinear problem}
Next, we return to the nonlinear problem.

\begin{assumption}\label{lastA}
Let $A:U\times(0,\infty)\to \mathcal M^{n\times n}_{{\rm sym}}$, $K:U\times(0,\infty)\to \mathcal M^{n\times n}$ and $B:U\times(0,\infty)\to \R^n$
be such that 
\begin{enumerate}[label=\tnum]
    \item $A(x,t)$ satisfies \eqref{Aelip2}, and 

    \item   $A(x,t)$ and $B(x,t)$ are bounded on $U\times (0,\infty)$.

    \item  There are constants $c_1\ge 0$ and $c_2\ge 0$ such that 
    \beq\label{condall}
- c_1|\xi|^2\le  \xi^{\rm T} K(x,t)\xi\le c_2|\xi|^2 \text{ for all $(x,t)\in U\times(0,\infty)$, all $\xi\in\R^n$.}
   \eeq 
\end{enumerate}
\end{assumption}

Let Assumption \ref{lastA} hold true for the rest of this section.
Condition \eqref{condall} in Assumption \ref{lastA} means that $K$ satisfies  \eqref{cond1} and \eqref{cond2} for all $T>0$. (In particular, if $K$ is bounded on $U\times(0,\infty)$, then \eqref{condall} certainly holds.)
By the boundedness of $B$ and $A$ on $U\times(0,\infty)$, there are positive numbers $M_0$ and $M_1$ such that 
\beq \label{Bmax}
|B(x,t)|\le M_0\text{ for all }(x,t)\in U\times(0,\infty), 
\eeq 
\beq \label{Tramax}
{\rm Tr}(A(x,t))\le M_1\text{ for all }(x,t)\in U\times(0,\infty).
\eeq 

We study the following initial boundary value problem
\beq\label{nonIBVP}
\begin{cases}
 \displaystyle   \frac{\partial u}{\partial t}-\langle A,D^2u\rangle +   u B\cdot \nabla u+ P'(u) (K\nabla u)\cdot \nabla u=0,&\text{ on }U\times(0,\infty),\\
    u(x,t)=u_*,&\text{ on } \Gamma\times(0,\infty),\\
    u(x,0)=u_0(x),&\text{ on }U,
\end{cases}
\eeq
where $u_*$ is a constant and $u_0(x)$ is a given function.

Define the nonlinear operator $L$ by \eqref{Ldef} for any function $u\in C_{x,t}^{2,1}(U\times(0,\infty))$ with the range of $u$ being a subset of $J$.

Assume $u\in C(\bar U\times[0,\infty))\cap C_{x,t}^{2,1}(U\times(0,\infty))$ is a solution of \eqref{nonIBVP} and satisfies 
\beq \label{ranu}
u(x,t) \in J\text{ for all }(x,t) \in U\times(0,\infty),
\eeq 

Same as  \eqref{reL}, we define  the linear operator $\mathcal L$ by 
\beqs %\label{reL}
\mathcal L w = w_t-\langle A(x,t),D^2w\rangle  +\widetilde B(x,t)\cdot \nabla w,\text{ where } \widetilde B(x,t)=u(x,t)B(x,t). 
\eeqs 
for any function 
$w\in C_{x,t}^{2,1}(U\times(0,\infty))$.

By the continuity of $u(x,t)$ on $\bar U\times [0,\infty)$  we must have
\beq \label{inib}
u(x,0)=u_* \text{ for }x\in \Gamma.
\eeq 
and, together with  the requirement \eqref{ranu}, $u_*\in \bar J$.
Also, the function $u_0(x)$, for $x\in U$, is continuous and $u(x,0)$  is its unique extension to a continuous function on $\bar U$. 
Hence, we can say $u(x,0)=u_0(x)$ on $\bar U$ and $u=u_*$ on $\Gamma\times[0,\infty)$.
Denote 
\beq
m_*=\min_{x\in \bar U} u(x,0)\text{ and } M_*=\max_{x\in \bar U} u(x,0).
\eeq
Then, thanks to \eqref{inib}, we have $m_*\le u_*\le M_*$. Since $u(U\times(0,\infty))\subset J$, it follows that  $m_*,M_*\in \bar J$, and, hence, the closed interval $[m_*,M_*]\subset \bar J$.

It follows the Maximum Principle -- Theorem \ref{maxprin} -- for all $T>0$ that 
\beq\label{umM}
m_*\le u(x,t)\le M_*\text{ on }\bar U\times [0,\infty).
\eeq

If $m_*=M_*$, then obviously 
\beq \label{triv}
u=m_*=u_*=M_*\text{ on } \bar U\times[0,\infty).
\eeq 
For this reason, we focus on the case $m_*<M_*$ now.
Pick any point $(x_0,t_0)\in U\times (0,\infty)$. Then $u(x_0,t_0)\in J\cap [m_*,M_*]$.

Consider the case $m_*\not \in J$. Since both $m_*,M_*$ are in the interval $\bar J$ and $m_*<M_*$, one deduces that $m_*$ cannot be a right-end point of $J$, hence $m_*$ must be the left-end point of $J$ and $\bar J$.
Similarly, if $M_*\not\in J$, then $M_*$ must be the right-end point of $J$ and $\bar J$.

With these discussions, we will refer to the following conditions, under the assumption $m_*<M_*$,  as
\begin{enumerate}[label={(E\theenumi)}]
\item\label{condc}
$M_*\in J$, $m_*\not\in J$ and  the right-hand limit $\displaystyle    \lim_{z\to m_*^+}P(z)$  exists and belongs to $\R\cup\{-\infty\}$,
\item\label{condd} 
$m_*\in J$, $M_*\not\in J$ and  the left-hand limit $\displaystyle    \lim_{z\to M_*^-}P(z)$ exists and belongs to $\R\cup\{\infty\}$.
\end{enumerate}

Consider  \ref{condc} and $\lambda> 0$. We can extend the function $e^{\lambda P(z)}$, for $z\in J$,  to a continuous function 
$E_\lambda:J\cup \{m_*\}\to [0,\infty)$.
For any function $F_\lambda$ in \eqref{ftrans},   we can extend it to a $C^1$-function, still denoted by $F_\lambda$, from $J\cup\{m_*\}$ to $\R$ by 
\beq\label{Fext1}
F_\lambda (s)=C\int_{s_0}^s E_\lambda(z) dz +C'\text{ for } s\in J\cup\{m_*\}.
\eeq

Consider  \ref{condd} and $\lambda< 0$. We similarly can extend any function $F_\lambda$ in \eqref{ftrans} to a a $C^1$-function, still denoted by $F_\lambda$ from $J\cup\{M_*\}$ to $\R$ by 
\beq\label{Fext2}
F_\lambda (s)=C\int_{s_0}^s E_\lambda(z) dz +C'\text{ for } s\in J\cup\{M_*\}.
\eeq

\begin{theorem}\label{Hthm2}
Suppose $m_*<M_*$. 
\begin{enumerate}[label=\tnum]
\item\label{interior} 
If $m_*,M_*\in J$, then there exist a number $C_0>0$ depending on $c_0$, $c_1$, $c_2$, $M_0$, $M_1$, $m_*$, $M_*$, and a number $\nu>0$ depending on $c_0$, $M_0$, $M_1$, $m_*$, $M_*$
 such that 
 \begin{align}\label{abs3}
    \max_{x\in \bar U} |u(x,t)-u_*| &\le C_0 e^{-\nu t }\max_{x\in \bar U}|u_0(x)-u_*|\quad \text{ for }t\ge 0.
    \end{align} 

\item\label{leftcase} 
Assume either $c_1=0$ 
or \ref{condc}. Then
      \begin{equation}\label{limu}
    \limsup_{t\to \infty} \max_{x\in \bar U} u(x,t)\le u_*.
    \end{equation}  
    If, in addition, $u_*=m_*$, then 
    \begin{equation}\label{stab}
    \lim_{t\to \infty} \max_{x\in \bar U} |u(x,t)-u_*| =0.
    \end{equation}  

\item\label{rightcase}
Assume either $c_2=0$ 
or \ref{condd}. Then 
    \begin{equation}\label{limi}
    \liminf_{t\to \infty} \min_{x\in \bar U} u(x,t)\ge u_*.
    \end{equation}  
    If, in addition, $u_*=M_*$, then one has \eqref{stab}.

    \item \label{mix}
    Consequently, if either
    \begin{enumerate}[label=\rnum]
        \item\label{aa} $c_1=0$ and \ref{condd}, 
     or
        \item\label{bb} $c_2=0$ and \ref{condc},
    \end{enumerate}
    then one has \eqref{stab}.
\end{enumerate}    
\end{theorem}
\begin{proof}
From \eqref{umM}, we have $|u(x,t)|\le \max\{|m_*|,|M_*|\}$ for all $(x,t)\in \bar U\times[0,\infty)$. 
Combining this with \eqref{Bmax} gives
\beq \label{tBmax}
|\widetilde B(x,t)|\le M_2 \text{ for all $(x,t)\in U\times(0,\infty)$, where 
$M_2= M_0\max\{|m_*|,|M_*|\}$.}
\eeq 

With $M_1$ and $M_2$ in \eqref{Tramax} and \eqref{tBmax}, let $s$, $\eta_*$ and $\nu$ be defined as in \eqref{sten}. Note that the last three numbers depend only on $c_0$, $M_0$, $M_1$, $m_*$, $M_*$.

Let 
\beq\label{llam}
\lambda_1>0\text{ and } \lambda_2<0\text{ such that }
\lambda_1\ge c_1/c_0 \text{ and } \lambda_2\le -c_2/c_0.
\eeq
For $j=1,2$, let $F_{\lambda_j}$ be defined by \eqref{ftrans} with $C=1$, $C'=0$ and $\lambda=\lambda_j$, and define 
\beq \label{wint}
w_j=F_{\lambda_j}(u)\text{ on }U\times(0,\infty).
\eeq 
 By Lemma \ref{q-lin}\ref{Fsub} and \ref{Fsuper}, we have 
 \beq \label{Lsign}
 \mathcal Lw_1\le 0\text{ and }\mathcal Lw_2\ge 0\text{  on }U\times (0,\infty).
 \eeq 

Below, whenever we apply Proposition \ref{Hthm1} to the operator $\mathcal L$, it means that $\widetilde L=\mathcal L$ with $b(x,t)=\widetilde B(x,t)$.

\medskip
\ref{interior} The proof of \eqref{abs3} is divided into three steps.

\medskip\noindent\emph{Step 1.} Since $m_*,M_*\in J$, we have 
\beq \label{bigrange}
u(\bar U\times[0,\infty))\subset [m_*,M_*]\subset J.
\eeq 
Select two numbers $\lambda_1$ and $\lambda_2$ that satisfy \eqref{llam}.
In this case, thanks to \eqref{bigrange}, we can extend to definition in \eqref{wint} to 
$w_j=F_{\lambda_1}(u)$ on $\bar U\times[0,\infty)$, for $j=1,2$.
Then we still have \eqref{Lsign}.

 For $j=1,2$, define $w_{*,j}=F_{\lambda_j}(u_*)$ and 
$\bar w_j=w_j-w_{*,j}$ on $\bar U\times[0,\infty)$.
Clearly,  $\bar w_j=0$ on $\Gamma\times[0,\infty)$ for $j=1,2$.
By Proposition \ref{Hthm1}(i) applied to  the operator  $\mathcal L$ and function $w:=\bar w_1$, we obtain from \eqref{upper2}, for $t\ge 0$, that
\beq \label{su1}
    \max_{x\in \bar U} \bar w_1(x,t) \le \eta_*^{-1} e^{-\nu t }\max\{0,\max_{\bar U}\bar w_1(x,0)\}
    \le \eta_*^{-1} e^{-\nu t }\max_{\bar U}|\bar w_1(x,0)|.
\eeq 
Similarly, by Proposition \ref{Hthm1}(ii) applied to the operator $\mathcal L$ and  function $w:=\bar w_2$, we have from \eqref{lower2}, for all $t\ge 0$, that
    \beq\label{su3}
    \max_{x\in \bar U} \bar w_2(x,t) \ge \eta_*^{-1} e^{-\nu t }\min\{0,\min_{\bar U}\bar w_2(x,0)\}
    \ge -\eta_*^{-1} e^{-\nu t }\max_{\bar U}|\bar w_2(x,0)|\}\text{ for }t\ge 0.
    \eeq 

\medskip\noindent\emph{Step 2.}     
The next step is to connect the inequalities \eqref{su1} and \eqref{su3} with $u(x,t)-u_*$. 
For that purpose, we denote 
$$C_1=\min\{e^{\lambda_1 P(m_*)},e^{\lambda_2 P(M_*)}\}\text{ and }
C_2=\max\{e^{\lambda_1 P(M_*)},e^{\lambda_2 P(m_*)}\}.$$
For $j=1,2$, we have 
$$0<C_1\le  F_{\lambda_j}'(z)=e^{\lambda_jP(z)}\le C_2 \text{ for }  z\in[m_*,M_*].$$
Above, we used the increasing property of $P$, see \eqref{Pd}.
Consequently, for $j=1,2$,
\beq \label{biLip}
C_1|s-u_*|\le |F_{\lambda_j}(s)-F_{\lambda_j}(u_*)|\le C_2|s-u_*| \text{ for }  z\in[m_*,M_*].
\eeq 
More specifically, by the Mean Value Theorem, one has,  for $j=1,2$,
\beq \label{suFsu}
\begin{aligned}
&C_1(s-u_*)\le F_{\lambda_j}(s)-F_{\lambda_j}(u_*)\le C_2(s-u_*) \text{ for } s\in[u_*,M_*],\\
&C_2(s-u_*)\le F_{\lambda_j}(s)-F_{\lambda_j}(u_*)\le C_1(s-u_*) \text{ for } s\in[m_*,u_*).
\end{aligned}
\eeq 
Consequently, when  $j=1$, we have from \eqref{suFsu},  for $s\in[m_*,M_*]$, that
\beq\label{uF1} s-u_*\le \max\{C_1^{-1}(F_{\lambda_1}(s)-F_{\lambda_1}(u_*)),C_2^{-1}(F_{\lambda_1}(s)-F_{\lambda_1}(u_*))\}.
\eeq 
When  $j=2$, we have,  for $s\in[m_*,M_*]$, 
\beq \label{uF2}
s-u_*\ge \min\{C_2^{-1}(F_{\lambda_2}(s)-F_{\lambda_2}(u_*)),C_1^{-1}(F_{\lambda_2}(s)-F_{\lambda_2}(u_*))\}.
\eeq

\medskip\noindent\emph{Step 3.} 
Now, combining inequality \eqref{uF1} with estimate \eqref{su1} yields, for any $t\ge 0$,
\begin{align*}
u(x,t)-u_*
&\le \max\{ C_1^{-1}(w_1(x,t)-w_{*,1}),C_2^{-1}(w_1(x,t)-w_{*,1})\}    \\
&\le \max \{  C_1^{-1}\eta_*^{-1} e^{-\nu t }\max_{\bar U} |w_1(x,0)-w_{*,1}|,C_2^{-1}\eta_*^{-1} e^{-\nu t }\max_{\bar U} |w_1(x,0)-w_{*,1}| \}\\
&\le   C_1^{-1}\eta_*^{-1} e^{-\nu t }\max_{\bar U}|w_1(x,0)-w_{*,1}|.
\end{align*}
Together with using \eqref{biLip} to estimate the last maximum, this yields
\beq\label{su2}
u(x,t)-u_*
\le  C_1^{-1} C_2\eta_*^{-1} e^{-\nu t }\max_{\bar U}|u(x,0)-u_*|.
\eeq 
Next, combining inequality \eqref{uF2} with estimate \eqref{su3} gives 
\begin{align*}
&u(x,t)-u_*
\ge \min\{ C_2^{-1}(w_2(x,t)-w_{*,2}),C_1^{-1}(w_2(x,t)-w_{*,2})\}    \\
&\ge \min \{ - C_2^{-1}\eta_*^{-1} e^{-\nu t }\max_{\bar U}|w_2(x,0)-w_{*,2}|,-C_1^{-1}\eta_*^{-1} e^{-\nu t }\max_{\bar U}|w_2(x,0)-w_{*,2}| \}\\
&\ge  -C_1^{-1}\eta_*^{-1} e^{-\nu t }\max_{\bar U}|w_2(x,0)-w_{*,2}|.
\end{align*}
Again, this and \eqref{biLip} imply
\beq\label{su4}
u(x,t)-u_*
\ge  -C_1^{-1}C_2\eta_*^{-1} e^{-\nu t }\max_{\bar U}|u(x,0)-u_*|.
\eeq 
Finally, combining estimates \eqref{su2} and \eqref{su4} yields
\beqs |u(x,t)-u_*|
\le  C_1^{-1}C_2 \eta_*^{-1} e^{-\nu t }\max_{\bar U}|u(x,0)-u_*|,
\eeqs 
which proves the desired estimate \eqref{abs3}.

\medskip
\ref{leftcase} We prove \eqref{limu} first.

\medskip\noindent\emph{Case 1: $c_1=0$.} In this case $\mathcal Lu\le Lu=0$. Then applying Proposition \ref{Hthm1}\ref{leftcase} to the operator $\mathcal L$ and function $w:=u-u_*$, we obtain \eqref{limu} from \eqref{basic1}.

\medskip\noindent\emph{Case 2: $c_1>0$ and \ref{condc}.}
We have the range $u(\bar U\times[0,\infty))$ is a subset of $J\cup  \{m_*\}$.
Let $\lambda_1=c_1/c_0>0$, and use the extended definition of the function $F_{\lambda_j}$ on $J\cup\{m_*\}$ given by \eqref{Fext1} with $C=1$, $C'=0$ and  $\lambda=\lambda_1$.
Then we can define  $w_{*,1}=F_{\lambda_1}(u_*)$ and 
$w_1=F_{\lambda_1}(u)$, $\bar w_1=w_j-w_{*,1}$ on $\bar U\times[0,\infty)$.

By \eqref{Lsign}, we have $\mathcal L\bar w_1\le 0$ on $U\times (0,\infty)$. 
By Proposition \ref{Hthm1}(i) applied to the operator $\mathcal L$ and  function $\bar w_1$, we have from \eqref{basic1} that
        \begin{equation}\label{stab1}
    \limsup_{t\to \infty} \max_{x\in \bar U} w_1(x,t) \le w_{*,1}.
    \end{equation}
By the increase and continuity of $F_{\lambda_1}$, it follows \eqref{stab1} that
      \begin{align*}
 F_{\lambda_1}(u_*)
 =w_{*,1}
 &\ge \limsup_{t\to \infty} \max_{x\in \bar U} F_{\lambda_1}(u(x,t))\\
&= \limsup_{t\to \infty} F_{\lambda_1}(\max_{x\in \bar U} u(x,t)) =F_{\lambda_1}(\limsup_{t\to \infty} \max_{x\in \bar U} u(x,t)) .
    \end{align*}
 Therefore, thanks to the strict increase of $F_{\lambda_1}$, we have
      \beq\label{ss1}
 u_*\ge \limsup_{t\to \infty}\max_{x\in \bar U} u(x,t),
    \eeq
which proves \eqref{limu}.  This completes the proof of \eqref{limu}.  

\medskip
Now, consider $u_*=m_*$. We have $u(x,t)\ge u_*$, then it follows \eqref{limu} that
$$ \limsup_{t\to \infty}\max_{x\in \bar U} |u(x,t)-u_*|=\limsup_{t\to \infty}\max_{x\in \bar U} (u(x,t)-u_*)\le 0.$$
Thus, we obtain \eqref{stab}.

\medskip
\ref{rightcase} We prove \eqref{limi} first.

\medskip\noindent\emph{Case 1: $c_2=0$.} In this case $\mathcal Lu\ge Lu=0$. Then \eqref{limi} follows \eqref{basic2} 
after applying Proposition \ref{Hthm1}\ref{rightcase} to the operator $\mathcal L$ and function $w:=u-u_*$.

\medskip\noindent\emph{Case 2: $c_2>0$ and \ref{condd}.}
The proof is the same as Part \ref{leftcase}, Case 2. Indeed, we have $u(\bar U\times[0,\infty))\subset J\cup \{M_*\}$.
Let $\lambda_2=-c_2/c_0<0$, and $F_{\lambda_2}$ be the extended function on $J\cup\{M_*\}$ defined by \eqref{Fext2} with $C=1$, $C'=0$ and $\lambda=\lambda_2$. Define $w_{*,2}=F_{\lambda_2}(u_*)$ and 
$w_2=F_{\lambda_2}(u)$, $\bar w_2=w_j-w_{*,2}$ on $\bar U\times[0,\infty)$.

By \eqref{Lsign},  we have $\mathcal L\bar w_2\ge 0$ on $U\times(0,\infty)$.
Then applying Proposition \ref{Hthm1}(ii)  to  the operator $\mathcal L$ and  function $\bar w_2$, we obtain from \eqref{basic2} that
\begin{equation}\label{stab2}
    \liminf_{t\to \infty} \min_{x\in \bar U}  w_2(x,t)\ge w_{*,2}.
    \end{equation}
Same as \eqref{ss1}, we have from \eqref{stab2} that
\beqs%\label{ss2}
  u_*\le   \liminf_{t\to \infty} \min_{x\in \bar U}u(x,t),
\eeqs 
which proves \eqref{limi}.

\medskip
Now that \eqref{limi} has already been established, consider $u_*=M_*$. We have  from \eqref{limi} that
$$ \limsup_{t\to \infty}\max_{x\in \bar U} |u(x,t)-u_*|
=-\liminf_{t\to \infty}\min_{x\in \bar U} (-|u(x,t)-u_*|)
=-\liminf_{t\to \infty}\min_{x\in \bar U} (u(x,t)-u_*)\le 0,$$
 hence, we obtain \eqref{stab} again.

 \medskip
 \ref{mix} On the one hand, one has, see  \eqref{absw},
 \beq\label{limabsu}
 \begin{aligned}
  &\limsup_{t\to \infty}\max_{x\in \bar U} |u(x,t)-u_*|\\
&  \le \max\{
 \limsup_{t\to \infty}\max_{x\in \bar U} (u(x,t)-u_*), 
 \limsup_{t\to \infty}\max_{x\in \bar U} (-(u(x,t)-u_*))
 \}\\
 &= \max\{
 \limsup_{t\to \infty}\max_{x\in \bar U} (u(x,t)-u_*), 
-\liminf_{t\to \infty}\min_{x\in \bar U} (u(x,t)-u_*) \}.
 \end{aligned}
 \eeq
On the other hand, in both cases \ref{aa} and \ref{bb}, one has \eqref{limu} and \eqref{limi}. 
Then combining \eqref{limabsu} with \eqref{limu} and \eqref{limi}, 
we obtain \eqref{stab}. 
 \end{proof}

As a consequence, we show below that the exponential rate of decay of $|u(x,t)-u_*|$, as $t\to\infty$, depends only on the asymptotic behaviors of $A(x,t)$ and $B(x,t)$ as $t\to\infty$, but not on the initial data $u_0(x)$ and the matrix $K(x,t)$.

\begin{corollary}\label{HCor}
Under Assumption \ref{lastA}, set 
\beq\label{lim1}
c_0'=\liminf_{t\to\infty} \inf_{x\in U}\min_{\xi\in\R^n,|\xi|=1} \xi^{\rm T} A(x,t)\xi ,
\eeq 
and let $M_0'$, $M_1'$ be two positive numbers such that 
\beq\label{lim3}
\limsup_{t\to\infty}\sup_{x\in U} |B(x,t)|<M_0'.
\eeq 
\beq\label{lim2}
\limsup_{t\to\infty}\sup_{x\in U} {\rm Tr}(A(x,t))<M_1',
\eeq 
If $m_*,M_*\in J$, then there is a number $\nu_*>0$ depending on $u_*$, $c_0'$, $M_0'$, $M_1'$ but not on the initial data $u_0(x)$,  such that 
    \begin{equation}\label{x2}
  \max_{x\in \bar U} |u(x,t)-u_*| =\mathcal O(e^{-\nu_* t}) \text{ as }t\to\infty.
    \end{equation} 
\end{corollary}
\begin{proof}
It is clear  from \eqref{abs3} that 
\beq\label{lim4}
\lim_{t\to\infty} \max_{x\in \bar U} |u(x,t)-u_*|=0.
\eeq
Note also that $u_*\in J$.
By \eqref{lim1}, \eqref{lim2}, \eqref{lim3}, \eqref{lim4}, there are $T>0$ and $m',M'\in J$ sufficiently close to $u_*$ with $m'<u_*<M'$ such that 
\beqs 
u(x,t)\in[m',M'] \text{ for all }(x,t)\in \bar U\times[T,\infty),
\eeqs  
\beqs
|B(x,t)|\le M_0',\quad 
{\rm Tr}(A(x,t))\le M_1'
\text{ for all }(x,t)\in U\times[T,\infty),
\eeqs
\beqs
\xi^{\rm T} A(x,t)\xi\ge \frac{c_0'}{2}|\xi|^2 \text{ for all }(x,t)\in U\times[T,\infty),\ \xi\in \R^n.
\eeqs
We repeat the proof of Theorem \ref{Hthm2}(i) with 
\beqs 
u:=u(x,t+T),\  A:=A(x,t+T),\  B:=B(x,t+T),\ K:=K(x,t+T),
\eeqs 
\beqs 
c_0:=c_0'/2,\ M_0:=M_0',\  M_1:=M_1',\ 
m_*:=m', \ M_*:=M'
\eeqs 
and the same numbers $c_1$, $c_2$. 
Note that the proof still works with the replacements $m_*:=m'$ and $M_*:=M'$ above even though $m'$ and $M'$ may not be the minimum and maximum  of $u(x,T)$ over $\bar U$.
Let $\nu_*=\nu$ be given by \ref{sten} with $M_2$ being replaced by $M_2'= M_0'\max\{|m'|,|M'|\}$, see \eqref{tBmax}. Then $\nu_*$ depends only on the numbers $c_0'/2$, $M_0'$, $M_1'$, $m'$, $M'$, hence is independent of $u_0(x)$. 
We obtain from \eqref{abs3} that 
\beqs 
|u(x,t+T)-u_*|
\le  C_* e^{-\nu_* t }\max_{\bar U}|u(x,T)-u_*|\text{ for some number } C_*>0.
\eeqs 
Hence, the estimate \eqref{x2} follows.
\end{proof}

\begin{example}\label{eg1}
Referring to Example \ref{fcase}, we consider cases \ref{Fe1} with \eqref{Fl1},  and \ref{Fe2} with \eqref{Fl20}, \eqref{Fl21} there. 
Either when  $J=[0,\infty)$ or $J=\R$, we always have $m_*,M_*\in J$. Therefore, for any $u_*\in J=\bar J$  and any corresponding solution $u$, one has from Theorem \ref{Hthm2}\ref{interior} the estimate \eqref{abs3} for all $t\ge 0$. 
\end{example}

\begin{example}\label{eg2}
Consider slightly compressible fluids as in case \ref{Fe3} of Example \ref{fcase}.  
We have  $J=(0,\infty)$, $u\ge 0$ on $\bar U\times[0,\infty)$, and $M_*\ge u_*\ge m_*\ge 0$.
Recalling \eqref{umM} and $u>0$ on $U\times(0,\infty)$, one has $M_*>0$, that is, $M_*\in J$.
Having \eqref{triv} in mind, we only consider $m_*<M_*$ below.

\medskip
\noindent\emph{Case 1: $u_*>0$.} We consider two subcases below.

Case 1a: $m_*>0$. Then $m_*,M_*\in J$ and one has from Theorem \ref{Hthm2}\ref{interior} the estimate \eqref{abs3} for all $t\ge 0$. 

Case 1b: $m_*=0$. Then (E1) is satisfied.
From Theorem \ref{Hthm2}\ref{leftcase}, we have the long-time estimate for $u$ by \eqref{limu} which is independent of $u_0(x)$.
If, in addition, $c_2=0$, then we have from Theorem \ref{Hthm2}\ref{mix}\ref{bb} the limit \eqref{stab}. 

\medskip
\noindent\emph{Case 2: $u_*=0$.} Then $m_*=0\not \in J$ and (E1) is satisfied. We have from Theorem \ref{Hthm2}\ref{leftcase} the limit \eqref{stab} which reads as
   \beq\label{eglim}
    \lim_{t\to \infty} \max_{x\in \bar U} u(x,t) =0.
    \eeq  
In fact, we can even obtain decaying estimates for all time.
Indeed, we can take $c_1>0$ and have from \eqref{Fl3} that $F_{\lambda_1}(s)=s^m$ with $m=c_1/c_0+1$.
In the proof of Theorem \ref{Hthm2}\ref{leftcase} at step \eqref{stab1}, we use the estimate \eqref{upper2} instead of the limit superior \eqref{basic1}. We obtain, for $t\ge 0$,
\beqs
    \max_{x\in \bar U} w_1(x,t)\le \eta_*^{-1} e^{-\nu t }\max_{x\in \bar U} w_1(x,0)
\eeqs
which yields
\beqs
    \max_{x\in \bar U} u^m(x,t)\le \eta_*^{-1} e^{-\nu t }\max_{x\in \bar U} u^m(x,0).
\eeqs
Therefore, instead of the limit \eqref{eglim}, we have 
$$  \max_{x\in \bar U} u(x,t)\le \eta_*^{-1/m} e^{-\nu t/m }\max_{x\in \bar U} u(x,0)\text{ for all }t\ge 0.$$
\end{example}

\medskip
\noindent\textbf{Data availability.} 
No new data were created or analyzed in this study.

\medskip
\noindent\textbf{Funding.} No funds were received for conducting this study. 

\medskip
\noindent\textbf{Conflict of interest.}
There are no conflicts of interests.

\bibliography{HI3Cite}{}
\bibliographystyle{plain}
\end{document}